\documentclass{article}
\usepackage[draft]{todonotes}   
\usepackage{amsthm,amsmath}
\usepackage{amscd,amssymb,latexsym,graphicx}
\usepackage[all]{xy}
\usepackage{cite}
\usepackage{mathtools}
\usepackage{hycolor}
\usepackage{xcolor}
\usepackage[
                       colorlinks=true,
                       linkcolor=black, %red,
                       citecolor=black, %green,
                       urlcolor=blue,
                     ]{hyperref}
\title{Classical Morse theory revisited -- I \\
         Backward $\lambda$-Lemma and
         homotopy type
        }
\author{Joa Weber\footnote{
        Financial support:
        FAPESP grant 2013/20912-4,
        FAEPEX grant 1135/2013, and 
        CNPq, Conselho Nacional de Desenvolvimento Cient\'{\i}fico
        e Tecnol\'ogico - Brasil.
  Instituto de Matem\'{a}tica, Estat\'{\i}stica
  e Computa\c{c}\~{a}o Scient\'{\i}fica,
  Universidade Estadual de Campinas,
  Rua S\'{e}rgio Buarque de Holanda~651,
  SP~13083-859 ,
  Campinas, Brasil.
        \hfill
        joa@ime.unicamp.br
        }
        \\
        IMECC UNICAMP }
\date{\today}
\newtheorem{theoremA}{Theorem}

\theoremstyle{definition}
\theoremstyle{remark}
\newcommand{\R}{{\mathbb{R}}}

\newcommand{\Dd}{{\mathcal{D}}}

\newcommand{\Gg}{{\mathcal{G}}}   % gauge transformations

\newcommand{\Ss}{{\mathcal{S}}}
\newcommand{\Tt}{{\mathcal{T}}}

\newcommand{\id}{{\rm id}}         % identity

\newcommand{\IND}{{\rm ind}}       % (Morse)index
 
\newcommand{\eps}{{\varepsilon}}

\def\NABLA#1{{\mathop{\nabla\kern-.5ex\lower1ex\hbox{$#1$}}}}
\def\Nabla#1{\nabla\kern-.5ex{}_{#1}}
\def\Tabla#1{\Tilde\nabla\kern-.5ex{}_{#1}}
\renewcommand{\Tilde}{\widetilde}

\newcommand{\p}{{\partial}}

\begin{document}
\maketitle
%%%%%%%%%%%% Abstract %%%%%%%%%%%%%%
\begin{abstract}
We introduce two tools, dynamical thickening and flow selectors,
to overcome the infamous discontinuity of the
gradient flow endpoint map near non-degenerate
critical points.
More precisely, we interpret the stable fibrations
of certain Conley pairs $(N,L)$, established
in~\cite{weber:2014c}, as a
\emph{dynamical thickening of the stable manifold}.
As a first application and to illustrate efficiency
of the concept we reprove a
fundamental theorem of classical Morse theory,
Milnor's homotopical cell attachment
theorem~\cite{milnor:1963a}.
Dynamical thickening leads to a
conceptually simple and short proof.
\end{abstract}%
%\tableofcontents

%\include{II}

%%%%%%%%%%%%%%%%%%%%%%%%%%%%%%%
%%%%%%%%%%%%%%%% Main Text  %%%%%%
%%%%%%%%%%%%%%%%%%%%%%%%%%%%%%%

Consider a connected smooth manifold
$M$ of finite dimension $n$. Suppose $f:M\to\R$
is a smooth function and $x$ is a non-degenerate
critical point of $f$ of Morse index $k$,
that is $df_x=0$ and
in local coordinates the Hessian matrix
$(\p^2 f/\p x^i\p x^j)_{i,j}$ at $x$ has precisely
$k$ negative eigenvalues, counting
multiplicities, and zero is not an eigenvalue.
Set $c:=f(x)$ and assume for simplicity
that the level set $\{f=c\}$ carries no
critical point other than $x$.

Morse theory studies how the topology
of sublevel sets $M^a=\{f\le a\}$ changes when
$a$ runs through a critical value $c$.
A fundamental tool is the concept of a flow,
also called a $1$-parameter group of
diffeomorphisms of $M$.
A common choice is the downward gradient
flow $\{\varphi_s\}_{s\in\R}$, namely the one
generated by the initial value problems
$\frac{d}{ds}\varphi_s=-(\nabla f)\circ\varphi_s$
with $\varphi_0=\id_M$. Existence is guaranteed,
for instance, if the vector field is of compact support.
Here $\nabla f$ denotes the gradient vector field
of $f$ on $M$. It is uniquely determined by
the identity $df(\cdot)=g(\nabla f,\cdot)$ after
fixing an auxiliary Riemannian metric $g$ on $M$.
Key properties of the downward gradient flow
are that $f$ decays along flow
lines $s\mapsto \varphi_s p$, for $p\in M$, and that
$\nabla f$ is orthogonal to level sets.
Consequently sublevel sets are forward
flow invariant. As $df_x=0$ $\Leftrightarrow$
$(\nabla f)_x=0$, any critical point $x$ is a
fixed point of the flow and non-degeneracy
translates into hyperbolicity.

By non-degeneracy of $x$ its unstable manifold
$W^u$ and descending disk $W^u_\eps$,
$$
     W^u=\{p\in M\mid\lim_{s\to-\infty}
     \varphi_sp=x\},\quad
    W^u_\eps=W^u\cap\{f\ge c-\eps\},
$$
are embedded open, respectively closed, disks in $M$ of
dimension $k=\IND(x)$; an embedding
$W^u_\eps\hookrightarrow M$ as a closed $k$-disk
exists only for every \emph{sufficiently small} $\eps>0$
(use the Morse-Lemma).
The boundary $S^u_\eps:=\p W^u_\eps$ is
called a descending sphere.
Consider instead the limit $s\to+\infty$
to get the stable manifold $W^s$
and ascending disk $W^s_\eps=W^s\cap\{f\le c+\eps\}$.
They have analogous properties
except that they are of codimension $k$.

In~\cite{weber:2014c}, see~\cite[Thm.~5.1]{Weber:2015c}
for details in the present finite dimensional case,
we implemented the structure of a disk bundle
on the compact neighborhood
$$
     N=N_x^{\eps,\tau}
     :=\left\{p\in M\mid\text{$f(p)\le c+\eps$,
     $f(\varphi_\tau p)\ge c-\eps$}\right\}
     _{\text{connected component of $x$}}
$$
of $x$ whenever $\eps>0$ is small and
$\tau>0$ is large.
The fibers are codimension-$k$ disks
with boundaries in the upper level set $\{f=c+\eps\}$
and parametrized by their unique point of
intersection, say $q^T$, with the unstable manifold.
The fiber over $x$ is $W^s_\eps$.
Each point of a fiber $N(q^T)$
reaches the lower level set
$\{f=c-\eps\}$ in time $T$ under the
downward gradient flow. Note that
$\{f=c-\eps\}$ intersects $W^u$
in the descending $(k-1)$-sphere
$S^u_\eps=\p W^u_\eps$.
Choose a tubular neighborhood $\Dd$ of
$S^u_\eps$ in $\{f=c-\eps\}$ to get a family
of codimension-$k$ disks $\Dd_q$, one for
each $q\in S^u_\eps$.
\begin{figure}%[b]
  \centering
  \includegraphics{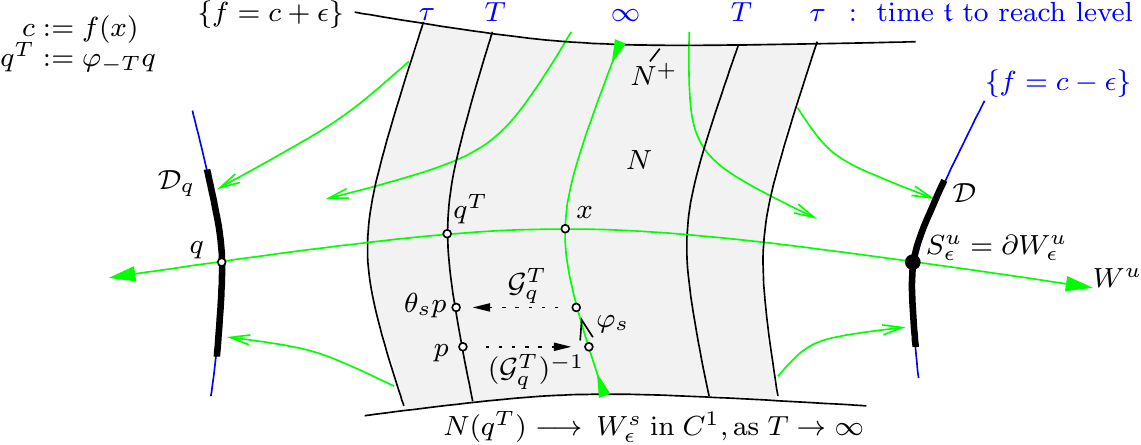}
  \caption{Dynamical thickening $(N,\theta)$
                 of the local stable manifold $(W^s_\eps,\varphi|)$
           }
  \label{fig:fig-N}
\end{figure}
By~\cite{weber:2014c,Weber:2015c}
we get a Lipschitz continuous ($C^{0,1}$) disk bundle
\begin{equation*}\label{eq:N_a-foliated}
     N
     =W^s_\eps
     \mathop{\dot{\cup}}_{T\ge\tau,q\in S^u_\eps}
      N(q^{T}),
     \quad
      N(q^{T})={\varphi_{T}}^{-1}(\Dd_q)
     \cap\{f\le c+\eps\},
\end{equation*}
over $\varphi_{-\tau} W^u_\eps$
which is $C^{1,1}$ away from the ascending disk
$W^s_\eps$. It is a key fact that the fibers
are diffeomorphic to $W^s_\eps$ via $C^1$
maps $\Gg^T_q:W^s_\eps\to N(q^T)$
which converge in $C^1$ to the identity on
$W^s_\eps$, as $T\to\infty$.
Furthermore, the fibration is forward flow invariant in the
sense that $\varphi_s$ maps a fiber $ N(q^T)$
into $ N(\varphi_sq^T)$.
Figure~\ref{fig:fig-N} illustrates the fibration and
the qualitative behavior of the forward
flow which is transverse to all fibers except the one
over $x$ which is invariant.
Conjugation by the diffeomorphism $\Gg^T_q$
provides on each fiber $N(q^T)$
a copy $\theta_s$ of the forward flow
$\varphi_s$ on $W^s_\eps$.
Now we reprove the cell attachment theorem.

\begin{theoremA}
[Milnor {\cite[I Thm.~3.2]{milnor:1963a}}]
Let $f:M\to\R$ be a smooth function, and let $x$
be a non-degenerate critical point with Morse
index $k$. Setting $f(x)=c$, suppose that
$f^{-1}[c-\eps,c+\eps]$ is compact and
contains no critical point of $f$ other than~$x$,
for some $\eps>0$. Then, for all sufficiently
small $\eps$, the set $M^{c+\eps}$ has the
homotopy type of $M^{c-\eps}$ with a
$k$-cell attached.
\end{theoremA}

\begin{proof}
Fix a Riemannian metric on $M$. 
Without loss of generality assume that $-\nabla f$
is of compact support,\footnote{
  Otherwise, substitute for $-\rho\nabla f$
  where $\rho:M\to\R$ is a smooth compactly
  supported cut-off function with $\rho\equiv 1$
  on the compact set $K:=f^{-1}[c-\eps,c+\eps]$.
  }
so it generates a flow
$\{\varphi_s\}_{s\in\R}$ on $M$.
Pick constants $\eps>0$ small and $\tau>0$ large
in order to meet the assumptions in~\cite{Weber:2015c} of
Theorem~5.4 (existence of the invariant fibration $N=N_x^{\eps,\tau}$)
and Definition~5.6 (induced fiberwise semi-flow $\theta$).
Figure~\ref{fig:fig-selector} illustrates the proof: First deform
$N\subset M^{c+\eps}$ along
$\theta$ towards the flow selector $\Ss^+$ and $W^u_\eps$,
then deform along $\varphi$.

\vspace{.1cm}
\textbf{0.~Definition of flow selector} (hypersurface
transverse to two flows): View
$$
     \Ss^+:=\{\varphi_{-\mathfrak{s}\circ\mathfrak{t}^-(p)}p\mid
     p\in\Ss^-\}\subset N
$$
as graph of a function $\mathfrak{s}\circ\mathfrak{t}^-$
over an open subset $\Ss^-\subset f^{-1}(c-\eps)$
where the coordinate lines are backward flow lines of $\varphi$
starting at $\Ss^-$ with coordinate the backward time.
By the flow box theorem this makes sense, as there
is no singularity of $\nabla f$ on $\Ss^-$. By the graph property
\emph{$\varphi$ will be transverse to $\Ss^+$}.

By~\cite[Thm.~1.2]{Weber:2015c} there is a $C^0$ time label function
$\mathfrak{t}:N\to[\tau,\infty]$, of class $C^1$ as a function
$N_\times:=N\setminus W^s\to[\tau,\infty)$, which assigns to
each point $p$ the time it takes to reach the lower level set
$f^{-1}(c-\eps)$ under the gradient flow $\varphi$. The hypersurface
$N^+:=\{p\in N\cap f^{-1}(c+\eps)\mid\mathfrak{t}(p)<\tau\}$
is called the \textbf{entrance set} of $N$ and
$N^+_\times:=N^+\setminus W^s$ its \textbf{regularization};
see Figure~\ref{fig:fig-N}. As each point of $N^+_\times$ hits
$f^{-1}(c-\eps)$ under $\varphi$ precisely once and transversely, the
corresponding subset $\Ss^-_\times\subset f^{-1}(c-\eps)$ is
diffeomorphic to $N^+_\times$. The \textbf{time label function}
$\mathfrak{t}^-:\Ss^-_\times\to(\tau,\infty)$ is defined by transfering
the time labels of $N^+_\times$. It is of class $C^1$.
Add the descending disk $S^u_\eps$ to define
$$
     \Ss^-:=\Ss^-_\times\mathop{\dot{\cup}} S^u_\eps
     =\{p\in f^{-1}(c-\eps)\mid N^+\cap\varphi_\R p\not=\emptyset\}
     \mathop{\dot{\cup}} S^u_\eps,\qquad
     \Ss^-_\times\stackrel{\varphi}{\cong} N^+_\times,
$$
as an open subset of $f^{-1}(c-\eps)$; see Figure~\ref{fig:fig-selector}.
Set $\mathfrak{t}^-=\infty$ on $S^u_\eps$.
\begin{figure}%[b]
  \centering
  \includegraphics{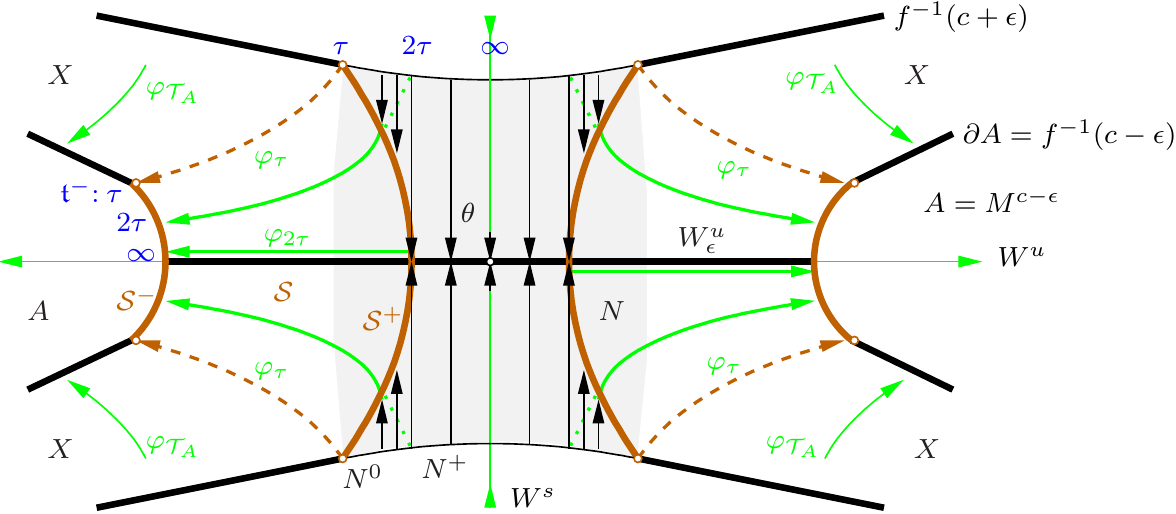}
  \caption{Flow selector $\Ss^+_\times=\Ss^+\setminus W^u$
          with transverse flows $\theta$ and $\varphi$
          }
  \label{fig:fig-selector}
\end{figure}
The function
$$
     \mathfrak{s}:(\tau,\infty)\to(\tau,2\tau),\quad
     \mathfrak{t}\mapsto 2\tau-\tau^2/\mathfrak{t},
$$
is smooth and extends continuously to  $[\tau,\infty]$ such that
$\mathfrak{s}(\tau)=\tau$ with $\mathfrak{s}^\prime(\tau)=1$
and $\mathfrak{s}(\infty)=2\tau$ with $\mathfrak{s}^\prime(\infty)=0$;
see Figure~\ref{fig:fig-selector} for the corresponding graph $\Ss^+$.
Observe that critical points of $\mathfrak{s}$
correspond precisely to tangencies of $\theta$ to the hypersurface
$\Ss^+_\times:=\Ss^+\setminus W^u$. But $\mathfrak{s}$ admits
no critical points on $(\tau,\infty)$, so \emph{$\theta$ is transverse to
$\Ss^+_\times$}. This proves that $\Ss^+_\times$ is a flow selector
with respect to $\varphi$~and~$\theta$.

\vspace{.1cm}
\textbf{I.~Strong deformation retraction} $r:M^{c+\eps}\to
M^{c-\eps}\cup W^u_\eps\cup X$ via $\theta$:
Let $\Ss$ be the region under the graph of $\Ss^+$,
that is the region bounded by $\Ss^-$ and $\Ss^+$
and the hypersurfaces indicated by dashed arrows
in Figure~\ref{fig:fig-selector}. The arrows are dashed to
indicate that they do not belong to $\Ss$, but to the closure $\bar\Ss$.
Consider the compact set
$X:=\left(f^{-1}[c-\eps,c+\eps]\setminus N\right)\cup\bar\Ss$ whose
boundary is given by $f^{-1}(c-\eps)$ and
$\left(f^{-1}(c+\eps)\setminus N^+\right)\cup\Ss^+$.
Deforming $N\setminus\bar\Ss$ along the flow lines of $\theta$
until the flow line hits either the flow selector $\Ss^+$ or the
descending disk $W^u_\eps$, while not moving the other points
of $M^{c+\eps}$ at all, defines the required strong deformation
retraction $r$. Continuity of $r$ holds since $\theta$ is transverse to
$\Ss^+_\times$ and $\Ss^+\setminus\Ss^+_\times=\Ss^+\cap W^u_\eps$
is reached under $\theta$ in infinite time just as is
$W^u_\eps\setminus \Ss$.

\vspace{.1cm}
\textbf{II.~Homotopy equivalence} $M^{c-\eps}\cup W^u_\eps\cup X
\sim M^{c-\eps}\cup W^u_\eps$ via $\varphi$:
Given the pair of closed sets $A:=M^{c-\eps}\subset \left(X\cup A\right)$,
consider the entrance time function $\Tt_A:X\cup A\to [0,\infty)$
which assigns to each point $p\in X\cup A$ the time it takes to reach
$A$ under $\varphi$. To see that $\Tt_A$ is well defined note that
$A$ and $X\cup A$ are both forward flow invariant under $\varphi$.
Indeed $\p A$ is a level set along which $-\nabla f$ is downward,
hence inward, pointing. The (topological) boundary of $X\cup A$ is
$\left(f^{-1}(c+\eps)\setminus N^+\right)\cup\Ss^+$
and $-\nabla f$ points inward along both pieces.

Given that $\varphi$ is transverse to $\p X$,
the function $\Tt_A$ is lower and upper
semi-continuous, hence continuous, because the
subset $A$ of $X\cup A$ is closed and
forward flow invariant, respectively;
cf.~\cite[Pf. of Thm.~B]{weber:2014c}.
Since $X$ is compact without critical points
$\Tt_A$ is bounded.  The map 
$h:[0,1]\times Z\to Z$ given by
\begin{equation*}
\begin{split}
     h(\lambda,p)=
     \begin{cases}
        p&
        \text{, $p\in A=M^{c-\eps}$,}
        \\
        \varphi_{\lambda\Tt_A(p)}p&
        \text{, $p\in X$,}
        \\
        \varphi_{\lambda4\tau^2/\mathfrak{t}(p)}p&
        \text{, $p\in\overline{W^u_\eps\setminus\Ss}=\varphi_{-2\tau} W^u_\eps$.}
     \end{cases}
\end{split}
\end{equation*}
is continuous as it is defined by three continuous parts which agree on overlaps:
$\Tt_A=0$ on $A\cap X$ and $\Tt_A=4\tau^2/\mathfrak{t}=2\tau$
on $\varphi_{-2\tau} S^u_\eps$.
The inclusion $\iota:A\cup W^u_\eps=:B\hookrightarrow Z:=X\cup A\cup
W^u_\eps$ and $h_1:=h(1,\cdot):Z\to B$ are reciprocal homotopy
inverses. Indeed $\iota\circ h_1=h_1\sim h_0=\id_Z$
and $h_1\circ\iota=h_1|_B\sim h_0|_B=\id_{B}$.
\end{proof}

% \footnote{
  Part two of $h_1$ unfortunately
  eliminates an outer piece of $W^u_\eps$ which we recover by
  $\varphi_{4\tau^2/\mathfrak{t}(\cdot)}(\cdot):\varphi_{-2\tau} W^u_\eps\to W^u_\eps$.
  So $h_1$ does not restrict to the identity on $W^u_\eps$,
  hence $h$ is not a deformation retraction
  of $X\cup A\cup W^u_\eps$ onto $A\cup W^u_\eps$.
%  }

%%%%%%%%%%%%%%%%%%%%%%%%%%%%%%%
%%%%%%%%%% Proof of theorem %%%%%%%%
%%%%%%%%%%%%%%%%%%%%%%%%%%%%%%%
%\subsubsection*{Proof}
%
%%% APPENDIX %%%%%%
%\appendix
%
% ACKNOWLEDGEMENTS %
% * various anonymous referees

%%%%%%%%%%%%%%%%%%%%%%%%%%%%%%%
\subsubsection*{Perspectives}
In the history of Morse theory
discontinuity of the flow trajectory end point
map $\varphi_\infty$ obstructed to carry out,
in a simple fashion, various constructions suggested
by geometry, for instance, to extend continuously
the inclusion map of an unstable manifold towards the closure.
It will be a future research project
to investigate the role of dynamical thickening
and flow selectors in such cases.

By~\cite{weber:2014c} dynamical thickening
can be defined in infinite dimensional contexts.

%%%%%%%%%%%%%%%%%%%%%%%%%%%%%%%
\subsubsection*{Added in proof}
Flow selector added
to correct the discontinuity in previous version.
Flow selectors arose in cooperation with 
Pietro Majer (2015) in two flavors - via Conley blocks
and via carving. Here we use a version of the
Conley block technique.

%%%%%%%%%%%%%%%%%%%%%%%%%%%%%%%%
%%%%%%%%% Acknowledgements %%%%%%%%
%%%%%%%%%%%%%%%%%%%%%%%%%%%%%%%%
\vspace{.1cm}
\noindent
{\small\bf Acknowledgements.} {\small
The author is grateful to %would like to thank
Stephan Weis for asking the right question just in time
and Kai Cieliebak for useful remarks concerning the flow selector.
}

%%%%%%%%%%%%%%%%%%%%%%%%%%%%%%%%%%%
%%%%%%%%%%%%%%%%%%%%%%%%%%%%%%%%%%%
%%%%%%%%%%%%%%%% References %%%%%%%%%%
%%%%%%%%%%%%%%%%%%%%%%%%%%%%%%%%%%%
%%%%%%%%%%%%%%%%%%%%%%%%%%%%%%%%%%%

%%%%%%%%%%%%%%%%%%%%%%%%%
%%%%%%%%% BibDesk %%%%%%%%%%
%%%%%%%%%%%%%%%%%%%%%%%%%
\bibliography{$HOME/Dropbox/0-Libraries+app-data/Bibdesk-BibFiles/library_math}{}

\begin{thebibliography}{1}

\bibitem{milnor:1963a}
J.~Milnor.
\newblock {\em Morse theory}.
\newblock Based on lecture notes by M. Spivak and R. Wells. Annals of
  Mathematics Studies, No. 51. Princeton University Press, Princeton, N.J.,
  1963.

\bibitem{weber:2014c}
J.~Weber.
\newblock {S}table foliations and semi-flow {M}orse homology.
\newblock {\em \href{http://arxiv.org/abs/1408.3842}{{\rm arXiv 1408.3842}},
  {\rm submitted}}, 2014.

\bibitem{Weber:2015c}
J.~Weber.
\newblock Contraction method and {L}ambda-{L}emma.
\newblock {\em S{\~a}o Paulo Journal of Mathematical Sciences}, 9(2):263--298,
  2015.

\end{thebibliography}
%$
%\bibliographysty le{plain}
         %   erzeugt:     [1] Joa Weber
\bibliographystyle{abbrv}
         %  erzeugt:      [1] J. Weber and 
%\bibliographystyle{alpha}
         %  article:    [Web05]  J. Weber
         %  book:      [Web05]  Joa Weber
         % more authors: [HZ87]

%%%%%%%%%%%%%%%%%%%%%%%%%
%%%%%%%%% standard %%%%%%%%%
%%%%%%%%%%%%%%%%%%%%%%%%%
%\begin{thebibliography}{00000}
%\small
%\end{thebibliography}

%\newpage
%%%%%%%%%%%%%%%%%%%%%%%%%%%%%%%%%%%%
%%%%%%%%%%%%% SYMBOLS %%%%%%%%%%%%%%
%%%%%%%%%%%%%%%%%%%%%%%%%%%%%%%%%%%%
%\subsection*{Symbols}
%$\N:=\{1,2,3,\dots\}$
%%%%%%%%%%%%%%%%%%%%%%%%%%%%%%%%%%%%
%%%%%%%%%% NOTIZEN %%%%%%%%%%%%%%%%%%%
%%%%%%%%%%%%% TODO's %%%%%%%%%%%%%%%%
%%%%%%%%%%%%%%%%%%%%%%%%%%%%%%%%%%%%
%\newpage
%\include{NOTIZEN}
%
%\listoftodos
%%%%%%%%%%%% To Do's %%%%%%%%%%%%%%%
%\vspace{1cm}\noindent
%{\Large\bf GENERAL TO DOs}
%\begin{itemize}
%\item 
%  MTC 
%\end{itemize}
%%%%%%%%%%%%%%%%%%%%%%%%%%
\end{document}